\documentclass[12pt]{article}

\usepackage{amssymb}
\usepackage{graphicx}
\usepackage[all]{xy}

\def \L {{\mathcal {L}}}

\def \d {{\mathrm {d}}}
\def \T {{\mathcal {T}}}

\def \b {{\bigskip}}
\def \m {{\medskip}}
\def \s {{\smallskip}}

\def \H {{\mathcal {H}}}
\def \L {{\mathcal {L}}}

\begin{document}

\author{Yvette Kosmann-Schwarzbach}

\title
{Beyond recursion operators}

\date{}

\maketitle

\noindent \textit{Abstract.} We briefly recall the history of 
the {Nijenhuis torsion} of $(1,1)$-tensors on manifolds and of the lesser-known {Haantjes torsion}. We then show how the {Haantjes manifolds} of Magri and the {symplectic-Haantjes structures} of Tempesta and Tondo generalize the classical approach to integrable systems in the bi-hamiltonian and symplectic-Nijenhuis formalisms, the sequence of powers of the recursion operator being replaced by a family of commuting Haantjes operators.

\smallskip

\begin{footnotesize}

\end{footnotesize}
\section{Introduction. Haantjes tensors generalize recursion operators}
The {Njenhuis torsion} of a $(1,1)$-tensor
 was defined in 1951 by Albert Nijenhuis, a student of the Dutch mathematician J. A. Schouten\footnote{A (1,1)-tensor on a manifold was called an 
{``affinor of valence two"}
 by Schouten and his contemporaries.
We find in the literature and we use indifferently the following expressions for a field of $(1,1)$-tensors: 
$(1,1)$-tensor [field], 
mixed tensor [of valence 2], 
field of endomorphisms of the tangent bundle, 
field of linear transformations, 
vector-valued [differential] $1$-form, 
$1$-form with values in the tangent bundle,
operator [on vector fields] [on $1$-forms].}, while the {Haantjes torsion} of a $(1,1)$-tensor was defined in 1955 by Johannes Haantjes, another of Schouten's students.
If the Nijenhuis torsion vanishes, the Haantjes torsion does also, but the converse is not true in general.
Since Nijenhuis tensors, i.e., $(1,1)$-tensors with vanishing Nijenhuis torsion, occur as {recursion operators} in the theory of integrable systems, 
one can expect the Haantjes tensors, i.e., $(1,1)$-tensors with vanishing Haantjes torsion, to play a role {``beyond recursion operators"}. 
\section{The search for differential concomitants and the Nijenhuis torsion}
\subsection{Schouten, Haantjes, Nijenhuis}
The mathematician {Jan A. Schouten} (1883-1971) is best known for his contributions to the modern form of the tensor calculus\footnote{See the article by his former and best known student, Nijenhuis \cite{N1972}.}, in particular for the differential concomitant that was later called ``the 
Schouten bracket" which he defined in an article that appeared in 
{\it Indagationes Mathematicae} in 1940 \cite{S1940}. His doctoral student,
{Johannes Haantjes} (1909-1956), defended his thesis at the University of    Leiden in 1933. Also a student of Schouten, {Albert Nijenhuis} (1926-2015) was awarded a doctorate at the University of Amsterdam in 1951. His article,
``$X_{n-1}$-forming sets of eigenvectors", appeared shortly thereafter in {\it Indagationes} \cite{N1951}, and four years later an article in two parts on   
``Jacobi-type identities for bilinear differential concomitants of certain tensor fields" was published in the same journal \cite{N1955}.
This second publication of Nijenhuis was preceded, only a few weeks earlier, by an article by Haantjes, ``On $X_m$-forming sets of eigenvectors", which also appeaed in {\it Indagationes} \cite{H1955}.

\subsection{The theory of invariants and the question of the integrability of eigenplanes} 

 The discovery of the Nijenhuis torsion followed a search for 
{differential concomitants} of tensorial quantities, which had its roots in the theory of invariants, going back to Cayley and Sylvester in the mid 19$^{\rm{th}}$ century, that led to
the Lie derivative in Sophus Lie's theory of continuous groups,
and later to the absolute differential calculus of Gregorio Ricci and Tullio Levi-Civita.
This search was extensively carried out by {Schouten} from the 1920's on. 
He wrote later \cite{S1954} that he had ``in 1940 succeeded in {generalizing Lie's operator} by forming a differential concomitant of two arbitrary contravariant quantities".
Then he disclosed the method he used to discover his concomitant: it was
by requiring that it be a derivation in each argument, which is the essential defining property of what is now called ``the Schouten bracket" of contravariant tensors.

Another field of enquiry was the search for conditions that ensure that, given a field of endomorphisms of the tangent bundle of a manifold,
assumed to have distinct eigenvalues, the {distributions spanned by pairs of eigenvectors} are integrable.

\subsection{The Nijenhuis torsion}
In 1951, Nijenhuis introduced a quantity 
defined by its components in local coordinates, 
$H^{ . .  \kappa}_{\mu \lambda}$, expressed
in terms of the components $h^{. \kappa}_\lambda$ of a $(1,1)$-tensor, $h$, and of their partial derivatives, 
 $${\it  H^{ . .  \kappa}_{\mu \lambda} =  2 
 h^{~. \rho}_{[ \mu}   \partial_{|\rho|}
 h^{. \kappa}_{\lambda]}  - 2 
  h^{. \kappa}_\rho \partial_{[ \mu}
   h^{. \rho}_{\lambda ]}}.$$
(The indices between square brackets are to be skew-symmetrized: the opposite term with these indices exchanged is to be added.) He then proved the tensorial character of this quantity \cite{N1951}. 
Because of the factor $2$, the $H^{ . .  \kappa}_{\mu \lambda}$ are actually the components of twice what is now called the  {\it  Nijenhuis torsion} of the $(1,1)$-tensor $h$, 
which is a skew-symmetric $(1,2)$-tensor, i.e., a vector-valued differential $2$-form.  

\m

\noindent {\it  Remark.} The name ``torsion" was adopted by Nijenhuis from the theory of 
complex manifolds, where the ``torsion" was defined for an almost complex structure by B. Eckmann and A. Fr\"olicher, also in 1951.
However, in the literature, the name ``Nijenhuis tensor" is often used for the ``Nijenhuis torsion".

\subsection{The Nijenhuis torsion without local coordinates}

It was also in his 1951 article that 
Nijenhuis introduced the {symmetric bilinear form}, depending on a pair of $(1,1)$-tensors, associated by polarization to the quadratic expression of the torsion. Then in 1955 \cite{N1955},
he introduced a {bracket notation} $[h,k]$ 
for this {symmetric bilinear form},
and
he found a {coordinate-independent} formula for this bracket. 
In particular the  {Nijenhuis torsion}, $\T_R = [R,R]$, of a $(1,1)$-tensor, $R$, on a manifold, $M$, is the {\it  $(1,2)$-tensor} $\T_R$ such that,
for all vector fields $X$ and $Y$ on $M$, 
$${\it  \T_R(X,Y) = [RX,RY]- R[RX,Y] - R[X,RY] + R^2[X,Y],}$$

\noindent
{\it  Remark.} We did not retain Nijenhuis's notation.
Our notation is simply related to his by $R = h$ and $\T_R = H$.

\subsection{The Fr\"olicher--Nijenhuis bracket}
In his 1955 article, Nijenhuis also defined what he called ``a concomitant for differential forms with values in the tangent bundle", that is, a graded bracket on the space of {vector-valued differential forms of all degrees}, extending the bilinear form associated to the torsion, and he proved that this bracket satisfies a graded Jacobi identity. (He also proved that  Schouten's brackets of contravariant tensors satisfy a graded Jacobi identity.)
This theory was soon developed in his joint article with A.~Fr\"olicher in 1956  \cite{FN1956}, and this graded Lie bracket became known as the {``Fr\"olicher--Nijenhuis bracket"}.

In a modern formulation, the Fr\"olicher--Nijenhuis bracket, $[U,V]_{FN}$, of a vector-valued $k$-form, $U$, and a
vector-valued $\ell$-form, $V$, is the vector-valued $(k+ \ell)$-form, 
$[U,V]_{FN}$, satisfying the equation 
$$
\L_{[U,V]_{FN}} = [\L_U, \L_V],
$$
where the bracket $[~,~]$ is the graded commutator of derivations of the algebra of differential forms, and where  $\L_U = [i_U,\d]$ and $\L_V = [i_V,\d]$ are the commutators of an interior product and the de Rham differential.
\m

Also in 1955, there appeared another, very different development of the theory of the Nijenhuis torsion of $(1,1)$-tensors, the {Haantjes torsion} of $(1,1)$-tensors.
\section{The Haantjes torsion}
\subsection{Haantjes (1909-1956)}
The Dutch mathematician {Johannes Haantjes}, 
after his doctoral defense in Leiden in 1933, was invited by  {Schouten} to join him as his assistant in Delft. 
From 1934 to 1938, they published several articles in collaboration, on spinors and their role in conformal geometry, and on the general theory of geometric objects, all in German except for one in English, papers that have 
almost never been cited.
After 1938, Haantjes was a lecturer at the Vrije Universiteit in Amsterdam. 
He was elected to the Royal Dutch Academy of Sciences in 1952, four years before his death at the age of 46.
He was among the ``distinguished European mathematicians" whom Kentaro Yano in 1982 recalled having met at the prestigious
International Conference on Differential Geometry organized in Italy in 1953  \cite{Y1982}. 
However, for half a century, extremely few citations of his work appeared in the literature, and the name of Haantjes was nearly forgotten.

\subsection{Haantjes's article of 1955}
In ``$X_m$-forming sets of eigenvectors" \cite{H1955}, Haantjes considered the case of a  field of endomorphims ``of class A", i.e., such that the eigenspace of an eigenvalue of multiplicity $r$ be of dimension $r$.
He introduced, in terms of local coordinates, a new quantity whose vanishing did not necessarily imply the vanishing of the Nijenhuis torsion
but was necessary and sufficient for the integrability of the distributions spanned by the eigenvectors.
From the Nijenhuis torsion $H$ of a $(1,1)$-tensor $h$ with components $H^{ . .  \kappa}_{\mu \lambda}$, he obtained the condition he sought as the vanishing of
$$ {\it  H^{ . .  \kappa}_{\nu \sigma} 
h^{\nu}_{. \mu} 
h^{\sigma}_{. \lambda} 
- 2 
H^{ . .  \sigma}_{\nu [ \lambda}
h^{\nu}_{. \mu ]} 
h^{\kappa}_{. \sigma} 
+ 
H^{ . .  \nu}_{\mu \lambda}
h^{\kappa}_{. \sigma} 
h^{\sigma}_{. \nu} 
.}
$$
These are the components of a $(1,2)$-tensor, twice the  {\it  Haantjes torsion} of the $(1,1)$-tensor $h$.
The components of the Haantjes torsion of $h$ are of degree $4$ in the components of $h$.

\subsection{First citations of Haantjes's article}
The 1955 article of Haantjes did not attract the attention of differential geometers or algebraists until the very end of the twentieth century. 
In fact, it was only cited twice before 1996!
In the twenty-first century, the ``Haantjes tensor" (i.e., in our terminology, the Haantjes torsion) started appearing, as an object of interest in algebra, in the work of Bogoyavlenskij \cite{B2004}, \cite{B2007}, and, mostly, in the {theory of integrable systems}. In 2007, in an article in {\it Mathematische Annalen}, Ferapontov and Marshall presented the Haantjes tensor as a ``differential-geometric approach to the integrability" of systems of differential equations, and 
 reformulated the main result of Haantjes's original paper as the theorem, 
 ``A system of hydrodynamic type with mutually distinct characteristic speeds is diagonalizable if and only if the corresponding Haantjes tensor [i.e., Haantjes torsion] vanishes identically" \cite{FM2007}.

\subsection{Haantjes torsion in coordinate-free form}
Changing notations, we denote a $(1,1)$-tensor by {\it  $R$}, its Nijenhuis torsion by $\T_R$, and we denote the Haantjes torsion of $R$ by  $\H_R$. We now formulate an intrinsic characterization of the {Haantjes torsion} of a $(1,1)$-tensor.

The {\it  Haantjes torsion} of a $(1,1)$-tensor $R$ is the {\it  $(1,2)$-tensor} $\H_R$ such that, for all vector fields $X$ and $Y$, 
$$\H_R(X, Y )= \T_R(RX,RY) - R(\T_R(RX,Y )) - R(\T_R(X,RY ))+ R^2(\T_R(X,Y )).$$
Explicitly,
$$\H_R(X,Y) = [R^2 X,R^2 Y] - 2 R  [R^2 X, RY] - 2 R  [RX, R^2 Y ]$$
$$ + 4 R^2 [RX,RY] +  R^2 [R^2 X,Y] +  R^2 [X, R^2 Y] - 2 R^3 [RX,Y]  - 2  R^3 [X,RY] + R^4 [X,Y].$$

Next, we shall generalize the definitions of the Nijenhuis torsion and of the Haantjes torsion of a $(1,1)$-tensor field on a manifold to any vector space equipped with a ``bracket". 
\section{Nijenhuis and Haantjes torsions associated to a ``bracket"}
\subsection{Definition} 
Let $\mu: E \times E \to E$ be a vector-valued skew-symmetric bilinear map on a real vector space $E$. 
For each linear map $R: E \to E$, 

\s

(i) the {\it  Nijenhuis torsion} of $R$  is the skew-symmetric $(1,2)$-tensor on $E$, denoted by $\T_R(\mu)$,  
such that, for all vectors $X$ and $Y$ in $E$, 

\m

\noindent {${\it  \T_R(\mu)(X,Y) = \mu(RX,RY)- R(\mu(RX,Y)) - R(\mu(X,RY)) + R^2(\mu(X,Y)),}$

\m

(ii)  the {\it  Haantjes torsion} of $R$ is the skew-symmetric $(1,2)$-tensor on $E$, denoted by $\H_R(\mu)$,
such that, for all vectors $X$ and $Y$ in $E$,
 
\m
 
\noindent $\H_R(\mu)(X,Y) =$

\s

\noindent{${\it  \T_R(\mu)(RX,RY)- R(\T_R(\mu)(RX,Y)) - R(\T_R(\mu)(X,RY)) + R^2(\T_R(\mu)(X,Y)).}$} 
\subsection{Lie algebroids}
The general definitions of the Nijenhuis torsion and of the Haantjes torsion are applicable when $E$ is the module of sections of a {\it  Lie algebroid}, $A \to M$, and 
$\mu$ is the Lie bracket of sections of $A$ (or to a
pre-Lie algebroid in which the bracket of sections does not necessarily satisfy the Jacobi identity), and $R$ is a section of $A \otimes A^*$. There are two important special cases:

(i) {\it  $E=TM$} is the module of vector fields on a manifold $M$, 
$\mu$ is the Lie bracket of vector fields and $R$ is a $(1,1)$-tensor, the case originally studied by Haantjes in 1955,

(ii) $E$ is a real {Lie algebra} and $R$ is a linear map.

\s

\noindent {\it  Remark.} In the case of $TM$ or, more generally, of a Lie algebroid $A$ over $M$, the Lie bracket of sections $\mu$ is {\it  only $\mathbb R$-linear, not $C^\infty(M)$-linear}. But the torsion $\T_R(\mu)$ of any $(1,1)$-tensor $R$ is $C^\infty(M)$-linear, i.e., it is a $(1,2)$-tensor.

\subsection{Haantjes torsion as torsion of the Nijenhuis torsion}

From the defining formula of the Haantjes torsion of a linear endomorphism $R$ of $E$ in terms of its Nijenhuis torsion we obtain immediately: 

\m

\noindent {\it Proposition.}
{The Haantjes torsion is related to the Nijenhuis torsion by}
$$\H_R(\mu) = \T_R(\T_R(\mu)).$$ 

This relation suggests the construction by iteration of higher Nijenhuis and Haantjes torsions of a linear endomorphism.

\subsection{Higher Nijenhuis torsions}Let $R$ be a linear endomorphism of a vector space $E$. Then $\T_R$ is the linear endomorphism of $E \otimes \wedge^2 E^*$ such that, for $\nu \in E \otimes \wedge^2 E^*$,
$$\T_R(\nu)  = \nu \circ (R \otimes R)   - R \circ \nu \circ (R \otimes {\mathrm{Id}})  - R \circ \nu \circ ( {\mathrm{Id}} \otimes R)  + R^2 \circ \nu.$$ 
For a vector space with bracket $\mu$, set $\T_R^{(1)}(\mu) = \T_R(\mu)$, which is, by definition, the Nijenhuis torsion $\T_R(\mu)$ of $R$.
Define
$$ \T_R^{(k + 1)}(\mu) =   \T_R(\T_R^{(k)}(\mu)), \, \, \, {\rm for} \, \, \, k \geq 1.$$
The $(1,2)$-tensors $\T_R^{(k)} (\mu)$ are of degree $2k$  in $R$.
We call the skew-symmetric $(1,2)$-tensors, $\T_R^{(k)} (\mu)$, 
for $k \geq 2$, the {\it   higher Nijenhuis torsions} of $R$.
For any skew-symmetric $(1,2)$-tensor $\mu$, and for all $k, \ell \geq 1$, 
$$
\T_R^{(k+\ell)}(\mu) = \T_R^{(k)}(\T_R^{(\ell)}(\mu)).$$

\subsection{Higher Haantjes torsions}
In the preceding notation,  
the Haantjes torsion of $R$ is 
$$\H_R(\mu) = \T_R(\T_R(\mu)) = \T_R^{(2)}(\mu).$$
Set $\H_R^{(1)}(\mu) = \H_R(\mu)$ and define 
$$ \H_R^{(k + 1)}(\mu) =   \T_R(\H_R^{(k)}(\mu)) , \, \, \, {\rm for} \, \, \, k \geq 1.$$
The  $(1,2)$-tensors $\H_R^{(k)} (\mu)$ are of degree $2(k+1)$ in $R$.
By definition, $\H_R^{(1)}(\mu)$ is the Haantjes torsion $\H_R(\mu)$ of $R$. 
We call the skew-symmetric $(1,2)$-tensors, $\H_R^{(k)} (\mu)$, 
for $k \geq 2$, the {\it  higher Haantjes torsions} which satisfy  
the very simple relation 
$$\H_R^{(k)} (\mu) = \T_R^{(k+1)} (\mu).$$
For any skew-symmetric $(1,2)$-tensor $\mu$, and for all $k, \ell \geq 1$, 
$$\H_R^{(k+\ell +1)}(\mu) = \H_R^{(k)}(\H_R^{(\ell)}(\mu)).$$
\subsection{A formula for the higher Haantjes torsions}
To a $(1,1)$-tensor,  Bogoyavlenskij associated a representation of the ring of real polynomials in 3 variables on the space of $(1,2)$-tensors \cite{B2004}.
Expanding the polynomial $(xy - zx - zy + z^2)^{k+1} = (z-x)^{k+1} (z-y)^{k+1}$ furnishes the general formula for the $(k+2)^2$ terms of the
{expansion of the $k$-th Haantjes torsion} of $R$, 
$$\H_R^{(k)}(\mu) (X,Y) =
\sum_{p=0}^{k+1}  \sum_{q=0}^{k+1}  (-1)^{2(k+1)-p-q} C_{k+1}^p C_{k+1}^q R^{p+q} \mu(R^{k+1-p}X,R^{k+1-q}Y).$$
To discover what roles, if any,  the higher Haantjes torsions can play in geometry and in the theory of integrable systems is an open question.
\subsection{Properties of the Nijenhuis and Haantjes torsions}
If a $(1,1)$-tensor field, $R$, on a manifold, $M$, is {diagonizable} in a local basis, 
$\displaystyle{(\frac{\partial}{\partial{x^i}})}$, $i = 1, \ldots, n$, 
with eigenvalues $\lambda_i(x^1, \ldots, x^n)$, its {Nijenhuis torsion} satisfies
$$\T_R (\mu)(\frac{\partial}{\partial{x^i}},\frac{\partial}{\partial{x^j}}) =  (\lambda_i  - \lambda_j) (\frac{\partial \lambda_i}{\partial{x^j}} \frac{\partial}{\partial{x^i}} + 
\frac{\partial \lambda_j}{\partial{x^i}} \frac{\partial}{\partial{x^j}}).$$  
Making use of the $C^{\infty}(M)$-bilinearity of $\T_R(\mu)$, it is easy to prove that, {if $R$ is diagonizable, the Haantjes torsion of $R$ vanishes}. 

\s
If there exists a basis of eigenvectors of $R$ at each point (in particular, if all the eigenvalues of $R$ are simple), the vanishing of the Haantjes torsion of $R$ is {\it  necessary and sufficient} for $R$ to be {\it  diagonalizable} in a system of coordinates.
\s

If {\it  $R^2 = \alpha \, {{\mathrm{Id}}}$}, where $\alpha$ is a constant, in particular, if $R$ is an {\it  almost complex structure}, i.e., when $R^2 = - {\mathrm{Id}}$,  then the Haantjes torsion 
is equal to the Nijenhuis torsion, up to a scalar factor, 

\s

\centerline{$\H_R(\mu) = 4 \alpha \T_R(\mu)$,}

\s
\noindent and, more generally,  $\H_R^{(k)}(\mu) =(4 \alpha)^k \T_R(\mu)$, for $k \geq 1$. 
\section{Haantjes manifolds and Magri--Lenard complexes}
\subsection{From Nijenhuis to Haantjes manifolds}
In a series of papers written since 2012, Franco Magri has defined the concept of a {Haantjes manifold}, demonstrated how the concept of a  Lenard complex on a manifold extends that of a Lenard chain associated with a bi-hamiltonian system}, related this theory to that of {Frobenius manifolds}, and developed applications to the study of differential systems \cite{M2012},  \cite{M2014},  \cite{M2015}, \cite{M2016}.

A {\it  Nijenhuis manifold} is a manifold endowed with a Nijenhuis tensor, i.e., a tensor whose Nijenhuis torsion vanishes.
In the theory of integrable systems, Nijenhuis tensors have also been called {\it Nijenhuis operators}, since they map vector fields to vector fields, as well as $1$-forms to $1$-forms, or {\it hereditary operators}, because they act as recursion operators. It is well known that every power of a Nijenhuis operator, $R$, is a Nijenhuis operator. Therefore, in a Nijenhuis manifold, the sequence of powers of the Nijenhuis tensor,
${\mathrm{Id}}, R, R^2 \ldots, R^k, \ldots$, is a {\it  family of commuting Nijenhuis operators}.

In the new framework, the role of this sequence of powers is played by a {\it family of Haantjes tensors}, i.e., $(1,1)$-tensors whose {Haantjes torsion vanishes}. Haantjes tensors are also called
{\it  recursion operators}.
A {\it  Haantjes manifold} is a manifold endowed with a {\it  commuting family  of Haantjes tensors}, $R_1, R_2, \ldots, R_k, \ldots$
In the examples, the family of Haantjes tensors is usually finite, 
in number equal to the dimension of the manifold, and {$R_1 = {\mathrm{Id}}$.}
We shall adopt the definition in this restricted sense.
\subsection{Magri--Lenard complexes}  
A {\it  Magri--Lenard complex} on a manifold, $M$, of dimension $n$, equipped with $n$ 
{\it  commuting $(1,1)$-tensors, $R_k$}, $k = 1,\ldots,n$, with $R_1 = {\mathrm{Id}}$, is defined by a pair $(X,\theta)$ such that

\m

(i) {the vector fields $R_k X$, $k = 1,\ldots,n$, \textit{commute pairwise}, 

\m

(ii) {the $1$-forms, $\theta R_k R_\ell$, $k, \ell = 1,\ldots,n$, are \textit{closed}}.

\m

\noindent In particular, $\theta$ itself is assumed to be closed and each 
$\theta R_k $ is a closed $1$-form.

\m
 
\noindent The notation  $\theta R$, where $\theta$ is a $1$-form and $R$ is a $(1,1)$-tensor, stands for $^t \! R \, \theta$, the $(1,1)$-tensor acting on the $1$-form by the dual map.

\m

Magri proved that, under a mild additional condition, if these properties are satisfied, {the operators $R_k$, $k = 1,\ldots,n$, are necessarily Haantjes tensors},
so that the underlying manifold of a Magri--Lenard complex is a {Haantjes manifold} \cite{M2017}.
\subsection{Magri--Lenard complexes generalize Lenard chains}  
In order to show how the Magri--Lenard complexes generalize the Lenard chains of bi-hamiltonian systems that were already defined by Magri in 1978 \cite{M1978}, we shall first recall how Nijenhuis operators appear in the theory of bi-hamiltonian systems. 

\m

If a vector field, $X$, leaves a  $(1,1)$-tensor, $R$, invariant, then,
$$0 = (\L_X R)(Y)= \L_X(RY) - R(\L_X Y) = [X,RY]  - R[X,Y],$$
for all vector fields $Y$. Therefore $R$, when applied to a symmetry $Y$ of the evolution equation, $u_t = X(u)$, yields a 
{new symmetry}, $RY$. 
If, in addition, $R$ is a 
{Nijenhuis operator}, 
applying the successive powers of $R$ yields a 
{sequence of commuting symmetries}, $R^k X$, $k \in \mathbb N$, known as a ``Lenard chain"\footnote{For the story of how the hierarchy of higher Korteweg--de Vries equations became known as a ``Lenard chain", named after Andrew Lenard (b. 1927), 
in papers by Martin Kruskal et al., see Lenard's letter reproduced in \cite{PS2005}.}
and therefore $R$ is a {\it  recursion operator} for each of the evolution equations in the hierarchy,
$u_t = (R^k X)(u)$.

\m

The geometric structure underlying the theory of integrable systems is the theory of Poisson--Nijenhuis manifolds, in particular the theory of symplectic--Nijenhuis manifolds. 
If $P_1$ and $P_2$ are 
{\it  compatible Hamiltonian operators}, i.e., Poisson bivectors such that their sum is a Poisson bivector, and if $P_1$ is invertible, i.e., defines a symplectic structure, then $R = P_2 \circ P_1^{-1}$ is a Nijenhuis operator. Thus, $(P_2, R)$ is called a ``Poisson--Nijenhuis structure" and $(P_1, R)$ is called a ``symplectic-Nijenhuis structure".
The theory of compatible Poisson structures originated in articles of  Gel'fand and Dorfman \cite{GD1980}, Fokas and Fuchssteiner \cite{FF1981}, Magri and Morosi \cite{MM1984}, and was further developed in \cite{KM1990} and \cite{KR2010}.

\m

We can now show that there is a Magri--Lenard complex associated to
a bi-hamiltonian system. Let $P_1$ and $P_2$ be compatible Hamiltonian operators. A vector field $X$ is called 
{\it  bi-hamiltonian} with respect to $P_1$ and $P_2$ if there exist 
{\it  exact} differential $1$-forms $\alpha_1 = \d H_1$ and $\alpha_2 = \d H_2$ such that 
$$X = P_1 (\alpha_1) = P_2 (\alpha_2).$$
Assume that $P_1$ is invertible, then the Nijenhuis operator $R = P_2 \circ P_1^{-1}$ generates a 
{sequence of commuting bi-hamiltonian vector fields}, $R^k X$, the so-called { Lenard chain}. 
The sequence of powers of $R$, $({\mathrm{Id}}, R, R^2 \ldots, R^k, \ldots)$, is a family of commuting Nijenhuis operators, and therefore a {\it  family of commuting Haantjes operators}. 
We set $\theta = \alpha_1$. Then $\theta R$ and all $\theta R^k$ are {closed} $1$-forms. Therefore, {\it the axioms of a Magri--Lenard complex are satisfied}.
In addition, the $1$-form $\theta$ and the recursion operator $R$ are invariant under~$X$.
\subsection{A Magri--Lenard complex on ${\mathbb R}^3$}
We present an example of a Magri-Lenard complex described by Magri in \cite{M2016}.
On ${\mathbb R}^3$ with coordinates $(u_1, u_2, u_3)$, consider the matrices 
$$K = \begin{pmatrix}{ 0 & 2 & 0 \cr - u_1 & 0 & 2 \cr - \frac{1}{2} u_2 & 0 & 0}\end{pmatrix} \, \, \, \, \rm{and} \, \, \, \, K^2 + u_1 {\mathrm{Id}} = \begin{pmatrix}{ - u_1 & 0 & 4 \cr - u_2 & - u_1 & 0 \cr
0 & - u_2 & u_1}\end{pmatrix}.$$
Matrices 
$K_0 = {\mathrm{Id}}, \, \, K_1 = K, \, \,   K_2 = K^2 + u_1 {\mathrm{Id}}$ commute.
\m 

Define
$\theta_0 = \theta = \d u_1$. We write $1$-forms as one-line matrices, and we consider the $1$-forms,

\s

$\theta_1 = \theta_{01} =  \theta K = 2 \d u_2$,                                                                                                                                                                                                                                                                                                                                                                                                                                                                                                                                                                                                                                                                                                                                                    
\, \, $\theta_2 = \theta_{02} = \theta K_2= - u_1 \d u_1 + 4 \d u_3$,

\s

$\theta_{11} = \theta_1 K = 2 \d u_1 $, \, \, \, 
$\theta_{12} = \theta_2 K = - 2 (u_2 \d u_1 + u_1 \d u_2)$,

\s

$\theta_{22} = \theta_2 K_2 = u_1^2 \d u_1 - 4 u_2 \d u_2$.

\m

\noindent All the $1$-forms, $\theta K_i K_j$, $0 \leq  i, j \leq 2$, are exact, and therefore {closed}.

\m

\noindent {\it  Remark.} Applying the successive powers of $K$ to $\theta$ {does not yield a sequence of closed forms.} While $\theta K^2$
and $\theta K^3$ are exact, $\theta K^4 = - 2 u_2\d u_1 - 4 u_1 \d u_2$ is not closed.

\b

Let $\displaystyle{X = \frac{\partial}{\partial{u_3}}}$.
The vector fields, $\displaystyle{X_0 = X =  \frac{\partial}{\partial{u_3}}, \, \, \, X_1 = K X = 2  \frac{\partial}{\partial{u_2}},}$ $\displaystyle{X_2 = K_2 X = 4  \frac{\partial}{\partial{u_1}} + u_1  \frac{\partial}{\partial{u_3}}}$, 
{commute}.

Therefore $\displaystyle{\left( {\mathbb R}^3, ({\mathrm{Id}}, K, K_2), \theta = \d u_1 , X =  \frac{\partial}{\partial{u_3}} \right)}$ is a Magri--Lenard complex.
In addition, $\L_X \theta = 0$ and $\L_X K = 0$.

\m

Computing the {Njenhuis torsion} $\T_K$ of $K$, we find that
$\T_K(e_ 1, e_2) =  e_2$, while
$\T_K(e_ 1, e_3) =  e_3$, and
$\T_K(e_ 2, e_3) = 0$. 
Thus, the vector-valued $2$-form $\T_K$ satisfies $\displaystyle{i_{\theta}\T_K = 0}$,
where $\theta = \d u_1$, but {$\T_K$ does not vanish}.
For a vector $X$ in ${\mathbb R}^3$, let $\T_K(X)$ be the endomorphism of ${\mathbb R}^3$ defined
by $Y \mapsto \T_K(X,Y)$. Then, 

\s

\noindent $\T_K(e_1) = \begin{pmatrix}{ 0 & 0 & 0 \cr0 & 1 & 0 \cr 0 & 0 & 1}\end{pmatrix}$, 
$\T_K(e_2) = - \begin{pmatrix}{ 0 & 0 & 0 \cr 1 & 0 & 0 \cr 0 & 0 & 0}\end{pmatrix}$, 
$\T_K(e_3) = - \begin{pmatrix}{ 0 & 0 & 0 \cr 0 & 0 & 0 \cr 1 & 0 & 0}\end{pmatrix}.$

\s

\noindent We now compute the Haantjes torsion of $K$.
Once we know the components of $\T_K$ and we have computed 
$K^2 = \begin{pmatrix}{ - 2 u_1 & 0 & 4 \cr - u_2 & - 2 u_1 & 0 \cr 0 & - u_2 & 0}\end{pmatrix}$,
we can compute the components of $\H_K$:
$$\H_K(e_1,e_2) = \T_K (Ke_1, Ke_2) - K \T_K(Ke_1, e_2) - K \T_K(e_1, Ke_2) + K^2 \T_K(e_1, e_2)$$
$$ = \T_K(- u_1 e_2 - \frac{1}{2} u_2 e_3, 2 e_1)- K \T_K(- u_1 e_2 - \frac{1}{2} u_2 e_3, e_2) - K \T_K(e_2, 2 e_1) - 2 u_1 e_2 - u_2 e_3$$
$$=  2 u_1 e_2 + u_2 e_3 - 2 u_1  e_2 - u_2 e_3 = 0,$$
and similarly, 
$\H_K(e_1,e_3) =  \H_K(e_2,e_3) = 0.$
Therefore $\H_K = 0.$

\s

Next, we compute the Nijenhuis and Haantjes torsions of $K_2$.
After computing 
$\displaystyle{(K_2)^2 = \begin{pmatrix}{(u_1)^2 & - 4 u_2 & 2 u_2 u_1 \cr 
0 & (u_1)^2  & - 4 u_2\cr
(u_2)^2  & 0 & (u_1)^2} \end{pmatrix}}$, 
we evaluate the {Nijenhuis torsion} of $K_2$ on the basis vectors and we obtain

\s

$\T_{K_2}(e_1, e_2) = u_2 e_3$,
$\T_{K_2}(e_1, e_3) = - 2 u_1 e_3$,
$\T_{K_2}(e_2, e_3) = 4 e_2$.

\s

\noindent Then we compute the {Haantjes torsion} of $K_2$ and we find that it vanishes. Therefore, $\displaystyle{\left( {\mathbb R}^3, ({\mathrm{Id}}, K, K_2)\right)}$ is a Haantjes manifold.

\m

{Why this example? }
The matrix $K$ in the preceding example is that of the {\it  integrable system of hydrodynamic type}, $U_t = K U_x$, where $U = \begin{pmatrix}{ u_1 \cr u_2 \cr u_3}\end{pmatrix}$, and $u_1$, $u_2$, $u_3$ are functions of two variables $(t,x)$. 
Explicitly, this differential system is   

\m

$\displaystyle{\frac{\partial u_1}{\partial t} = \frac{\partial u_2}{\partial x},}$

\s

$\displaystyle{\frac{\partial u_2}{\partial t} = -  u_1 \frac{\partial u_1}{\partial x}   + 2  \frac{\partial u_3}{\partial x},}$
 
 \s
 
$\displaystyle{\frac{\partial u_3}{\partial t} =  -  \frac{1}{2} u_2 \frac{\partial u_1}{\partial x}.}$

\m

Another case to which the geometric structure of Haantjes manifolds is applicable is that of the {\it  dispersionless Gelfand--Dickey equations} defined by the $(1,1)$-tensor, 
$K = \begin{pmatrix}{ 0  & 1 & 0 \cr u_1 & 0 & 1\cr u_2  & u_1 & 0 }\end{pmatrix}.$
\section{WDVV equations and Magri--Lenard complexes}
Magri showed how the geometric structures on Haantjes manifolds are related to the solutions of the {\it  WDVV equations}\footnote{This system of partial differential equations is named after E. Witten, R. Dijkgraaf, E. Verlinde and H. Verlinde.}. These are the equations satisfied by the partial derivatives of the Hessian, i.e., the matrix of second-order partial derivatives,
of a function $F$ of $n$ variables, $(x^1, x^2, \ldots, x^n)$. 
Let the {Hessian matrix} of $F$ be denoted by $h$ and assume that the matrix $\displaystyle{\frac{\partial h}{\partial x^1}}$ 
is invertible. The WDVV equations can be written as the set of nonlinear equations,
$$ \frac{\partial h}{\partial x^i}\left ( \frac{\partial h}{\partial x^1}\right )^{-1}\frac{\partial h}{\partial x^j} = \frac{\partial h}{\partial x^j} \left ( \frac{\partial h}{\partial x^1} \right )^{-1}\frac{\partial h}{\partial x^i}, \, \, \, \hspace{.2cm} i, j = 1, \ldots, n.$$
They express the pairwise {commutativity} of the matrices 
$$\left ( \frac{\partial h}{\partial x^1}\right )^{-1} \frac{\partial h}{\partial x^i} , \, \, \, \hspace{.2cm} i = 1, \ldots, n.$$
 
\b

Given a solution, $F$, of the WDVV equations, consider the $1$-forms
$\theta_{ij} = \d a_{ij}$, $i, j = 1 \ldots, n$, 
where the $\displaystyle{a_{ij} =  \frac{\partial^2F}{\partial x^i \partial x^j}}$  are the entries of the Hessian matrix $h$ of $F$. 
Assume that the $1$-forms $\theta_{1j}$, $ j = 1, \ldots, n$, are linearly independent, and define operators $R_k$ by the condition 
$$\theta_{1j} R_k = \theta_{jk}.$$
Then $R_1= {\mathrm{Id}}$ and
$$\theta_{11} R_i R_j = \theta_{1i} R_j = \theta_{ij}.$$
\noindent{\it  Proposition.} Consider the {commuting vector fields} $\displaystyle{X_k=  \frac{\partial}{\partial x^k}}$. 
Then the operators $R_k$ satisfy the relation 
$$X_k= R_k  \frac{\partial}{\partial x^1}.$$
\noindent{\it Proof.} On each of the linearly independent $1$-forms $\theta_{1j} = \d a_{1j}$, $j = \ldots, n$,
the vector fields $\displaystyle{X_k = \frac{\partial}{\partial x^k}}$ 
and $\displaystyle{R_k  \frac{\partial}{\partial x^1}}$ take the same value, 
$\displaystyle{\frac{\partial a_{1j}}{\partial x^k} =\frac{\partial a_{jk}}{\partial x^1}}$.

\b

The operators $R_k$ commute because $F$ is assumed to be a solution of the WDVV equations. In fact, $\displaystyle{ R_k \frac{\partial h}{\partial x^1} = \frac{\partial h}{\partial x^k}}$. 
Therefore the operators $R_k$, the vector field $\displaystyle{\frac{\partial}{\partial x^1}}$ and the $1$-form $\theta_{11}$ define a {\it  Magri--Lenard complex}. 

\m

Conversely, consider a {\it  Magri--Lenard complex} $(M, R_k, X, \theta)$.
{Locally}, on an open set of the manifold $M$, the commuting vector fields 
$X_k = R_k X$ define {coordinates} $x^k$, and the closed $1$-forms $\theta_{ij} = \theta R_i R_j$ admit local {potentials} $a_{ij}$,
$$\theta_{ij} =  \d a_{ij}.$$
For $i, j, k = 1, \ldots , n$, consider the functions 

\centerline{$c_{ijk} = \langle \theta_{ij}, X_k \rangle = \langle \theta R_i R j, R_k X \rangle  =  \langle \theta, R_i R j R_k X \rangle.$}

\m

\noindent In local coordinates,

\centerline{$\displaystyle{c_{ijk} = \langle \theta_{ij}, X_k \rangle = \langle \d a_{ij}, \frac{\partial}{\partial x^k} \rangle = \frac{\partial a_{ij}}{\partial x^k}}.$}

\m

\noindent Because the operators $R_k$ commute pairwise, functions 
$c_{ijk}$ are symmetric. Therefore the functions $\displaystyle{\frac{\partial a_{ij}}{\partial x^k}}$ are {symmetric}, which implies that the  $a_{ij}$ are the second-order partial derivatives 
of a function $F(x^1, \ldots, x^n)$, $\displaystyle{a_{ij} =  \frac{\partial^2 F}{\partial x^i \partial x^j}.}$ 
Then the Hessian of $F$ {satisfies the WDVV equations}.
\section{Lenard--Haantjes chains}
To conclude this survey of modern work based on the 1955 article of Haantjes, I must mention recent work of 
Tempesta and Tondo \cite{TT2016a}, \cite{TT2016b}, \cite{TT2017}.

\m
 
A {\it  symplectic-Haantjes manifold} of dimension $2n$ is a symplectic manifold $(M, \omega)$ endowed with a 
family 
$(K_0 = \mathrm{Id}, K_1, \ldots, K_{n-1})$ of $n$ linearly independent Haantjes tensors such that:

(i) each map $\omega^\flat \circ K_i: TM \to T^*M$, $i = 0, \ldots, n - 1$, 
is skew-symmetric,

(ii) the $K_i$'s, $i = 0, \ldots, n - 1$, generate a $C^\infty(M)$-module of Haantjes tensors,

(iii) for all $i, j = 0, \ldots, n - 1$, $K_i K_j$ has a vanishing Haantjes torsion, and  $K_i K_j = K_j K_i$.

\m

``Lenard--Haantjes chains" are constructed from a given Hamiltonian $H$ on a symplectic-Haantjes manifold by defining Hamiltonians $H_j$ such that
$$\d H_{j+1} = \d H \,  K_j.$$
Then the Poisson bracket of any two Hamiltonians $H_j$ in the chain vanishes. 

\m

Among the examples given by Tempesta and Tondo are the generalized St\"ackel systems, where $\omega$ is the canonical symplectic form on $T^*({\mathbb R}^{n})$ with coordinates $(q^i, p_i)$, 
and the $K_j$'s are diagonal operators defined in terms of the cofactors of a St\"ackel matrix, an invertible matrix whose $i$-th row depends only on the coordinate $q^i$.

\m

Many other applications of the Haantjes tensors can be found in the publications and in the preprints of Tempesta and Tondo, as well as in Magri's articles, published or in progress. The comparison of the methods thus proposed to investigate the geometry of integrable systems remains to be done.

\b

\noindent {\it Acknowledgments.} I would like to thank the organizers of the 37th Workshop on Geometric Methods in Physics for their invitation and for an interesting meeting. I must thank Giorgio Tondo for communicating  to me unpublished material on his joint work with Piergiulio Tempesta, and I extend my sincere and friendly thanks to Franco Magri upon whose work, both published and in preprint form, I have largely drawn in this paper.

\begin{small}

\end{small}

\bigskip

\noindent Yvette Kosmann-Schwarzbach

\noindent Paris, France

\noindent e-mail: {\tt{yks@math.cnrs.fr}}
\end{document}